\documentclass{article}
\usepackage{amsmath}
\usepackage{ntheorem-hyper}
\usepackage{amssymb}
\usepackage{authblk} 
\usepackage{algorithm}
\usepackage{algorithmicx}
\usepackage{algpseudocode}
\usepackage{graphicx}
\usepackage{graphbox}
\usepackage{wrapfig}
\usepackage{subfigure}
\usepackage{svg}
\usepackage{epstopdf}
\usepackage{tabularx}  
\usepackage{booktabs}
\usepackage[numbers, sort]{natbib}
\usepackage{hyperref}
\usepackage{xcolor} 
\usepackage{color}
\usepackage[title]{appendix}
\usepackage{longtable}

\usepackage[T1,T2A]{fontenc}
\usepackage[utf8]{inputenc}
\usepackage[russian, german, english]{babel}

\hypersetup{urlcolor=blue, citecolor=red}

\def\RR{\mathbb{R}}

\title{Supersqueezed sphere}
\author[+]{Hongda Qiu}

\begin{document}

\maketitle

\section{Introduction}

Let $S$ be a smooth surface in $\RR^3$ with normal curvatures at most 1 in absolute value.
In other words, $S$ is
not allowed to bend more than the unit sphere.

\smallskip

\textit{Assume $S$ is a topological sphere.
Can the volume enclosed by $S$ be less than the volume of the unit ball?}

\smallskip

This question was posed several years ago by Dmitri Burago and Anton Petrunin
\cite[p.~11]{petrunin2023pigtikalpuzzlesgeometryi}\cite{anton2022mathoverflow}.

Soon we will construct such an example,
but before getting there, let us talk a bit about the history of the problem.

Consider the analogous question in the plane:
\smallskip

\textit{Can a closed simple plane curve with curvature at most 1 in  absolute value bound a smaller area than the unit disc?}

\begin{wrapfigure}{r}{36 mm}
\vskip-8mm
\centering
\includegraphics[width=32mm]{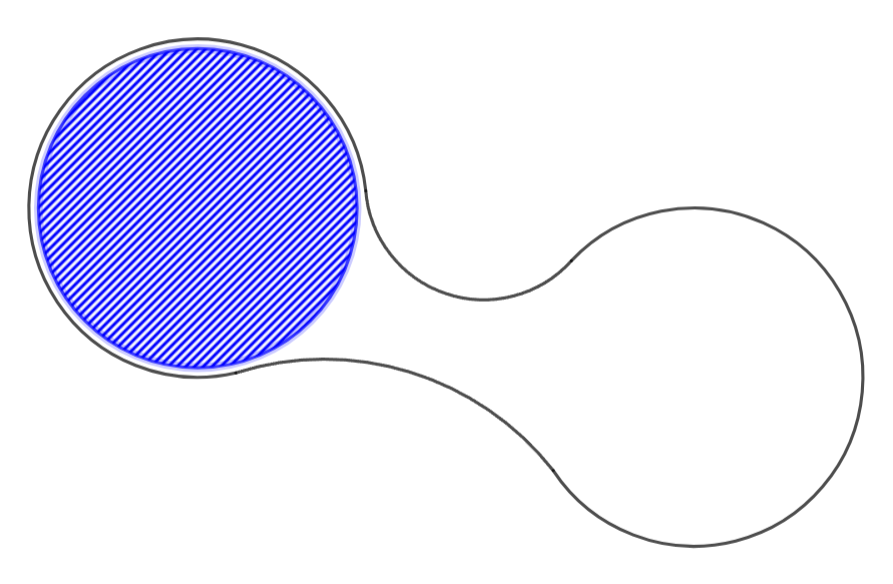}
\vskip0mm
\end{wrapfigure}

\smallskip

Such curves do not exist.
In fact, German Pestov and Vladimir Ionin \cite{pestov1959} proved that such a curve necessarily encloses a unit disc,
which is of course a stronger statement.
(Note that we do not assume convexity.)

A quick corollary of this theorem is that a surface of revolution with normal curvatures $\leq 1$ in absolute value necessarily encloses a unit ball in $\RR^3$, implying a partial answer to our main question.
There are several other partial answers.
For instance, the case when the surface is contained in a ball of radius $2$ follows from my result in \cite{hongda2025ball}.
In particular, one cannot obtain an example from small variations of the unit sphere.
There are also no counterexamples among star-shaped surfaces; this was proved by Terence Tao in collaboration with AI \cite{anton2022mathoverflow}.

All this leads to the following stronger version of our question.

\smallskip

\textit{Does $S$ enclose a unit ball?}

\smallskip

This question was considered in several papers by Vladimir Lagunov and his advisor Abram Fet \cite{lagunovfet1961russianI}\cite{lagunovfet1961russianII}\cite{lagunov1959russian}\cite{lagunov1960russianI}\cite{lagunov1961russianII}.\footnote{Partial translations of these papers were made by Richard L. Bishop \cite{lagunov1963extremalproblemssurfacesprescribed}.}
They constructed counterexamples and showed that $S$ must enclose a ball of radius $r=\frac{2\sqrt{3}}{3}-1$;
this is an optimal bound for all closed surfaces with normal curvatures at most 1 (not necessarily spheres).
For spheres they proved that the optimal bound is $r=\frac{\sqrt{6}}{2}-1$.
{

\begin{wrapfigure}{r}{16 mm}
\vskip-6mm
\centering
\includegraphics[width=15 mm]{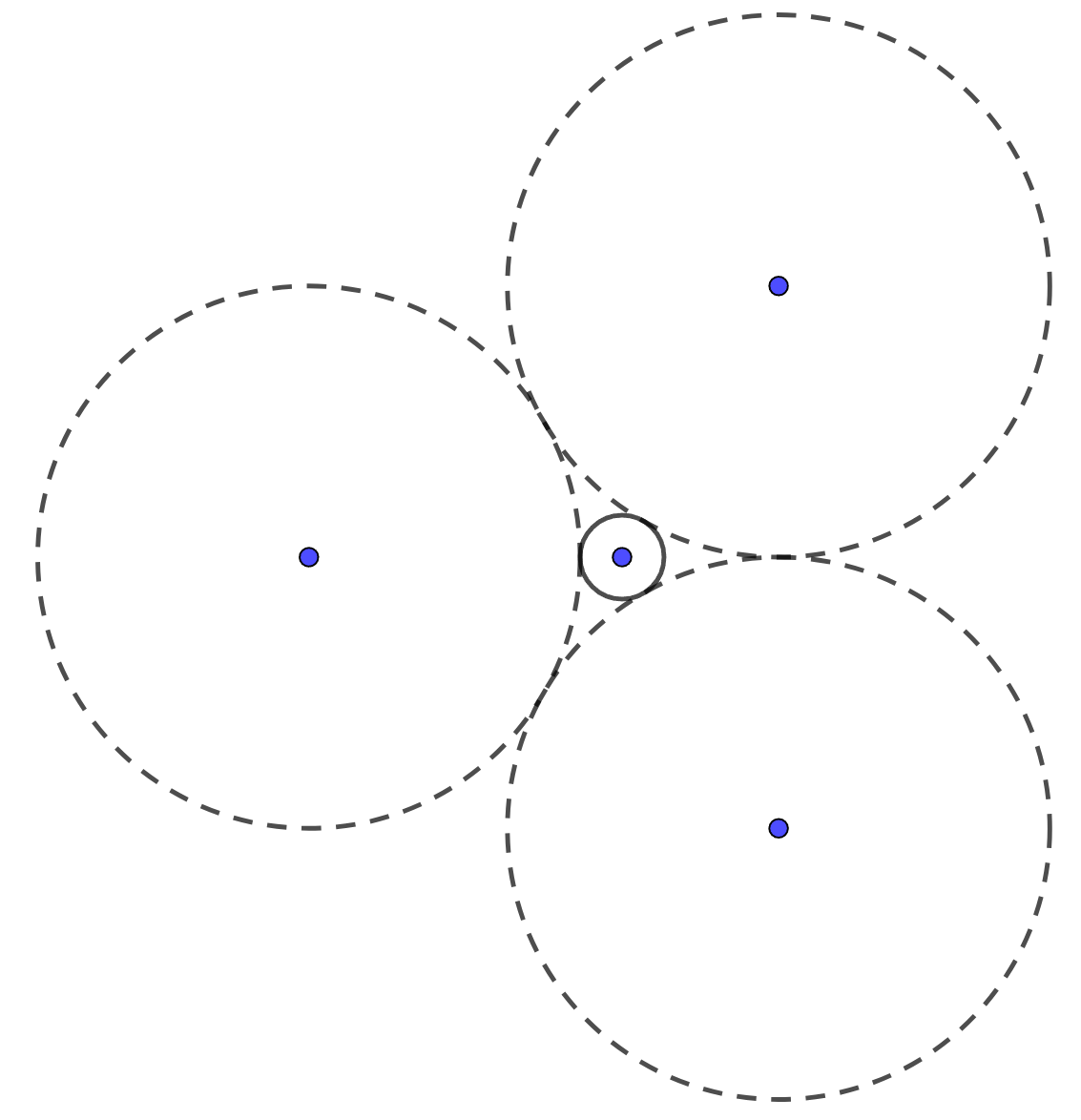}
\vskip0mm
\end{wrapfigure}

Geometrically, the former number is the radius of the maximal disc inscribed in three unit circles as shown in the figure,
and the latter is the radius of the maximal ball squeezed between four mutually tangent unit spheres.

}

\begin{wrapfigure}{r}{76 mm}
\vskip-6mm
\centering
\includegraphics[width=75 mm]{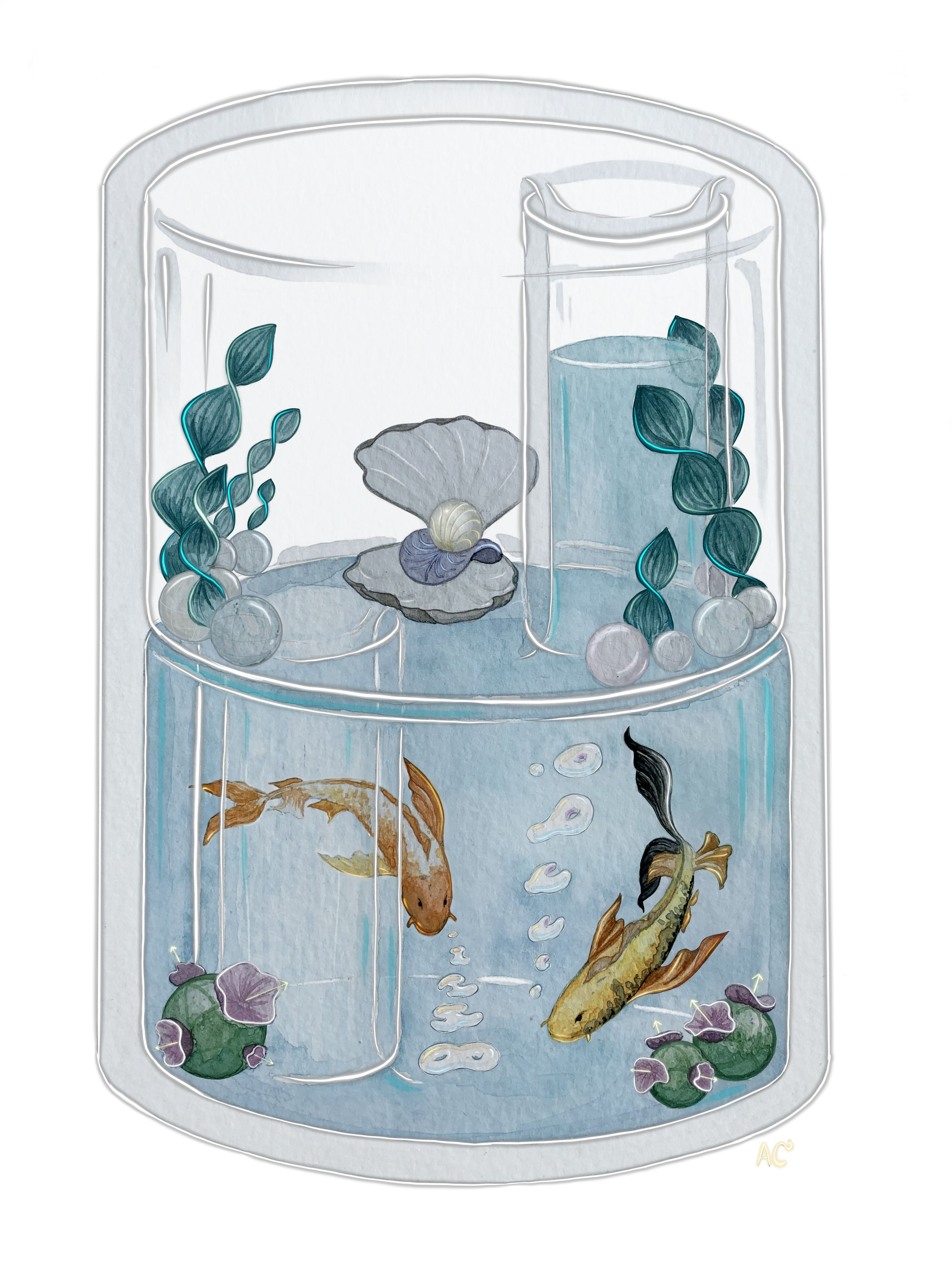}
\vskip0mm
\end{wrapfigure}

These examples are often referred to as \textit{Lagunov's fishbowls}; the drawing by Ana Cristina Chávez Cáliz explains why. 
Detailed discussions of this example are given by Anton Petrunin in his textbook \cite[Section~11.C]{petrunin2024differentialgeometrycurvessurfaces} and by Sergio Zamora Barrera in his video lesson \cite{Barrera2023}. You may look at the drawing and realize that the enclosed body can be made arbitrarily thin except for the belt where three pieces of the surface come close together --- this is the place where volume lives.

The surface shown has genus 2 (it is a sphere with two handles),
but a similar example was constructed for the sphere.
The resulting surface can be thought of as a fattened Bing's house \cite{bing1964}. It seems that Lagunov and Fet did not know this example, and had to reinvent it. We are skipping the details since we will have to use the same trick in our construction.

\section{Jerrycan}

\begin{figure}[ht!]
\centering
\subfigure{
\includegraphics[width=0.40\linewidth]{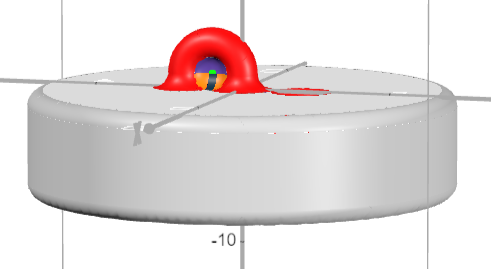}
}
\subfigure{
\includegraphics[width=0.45\linewidth]{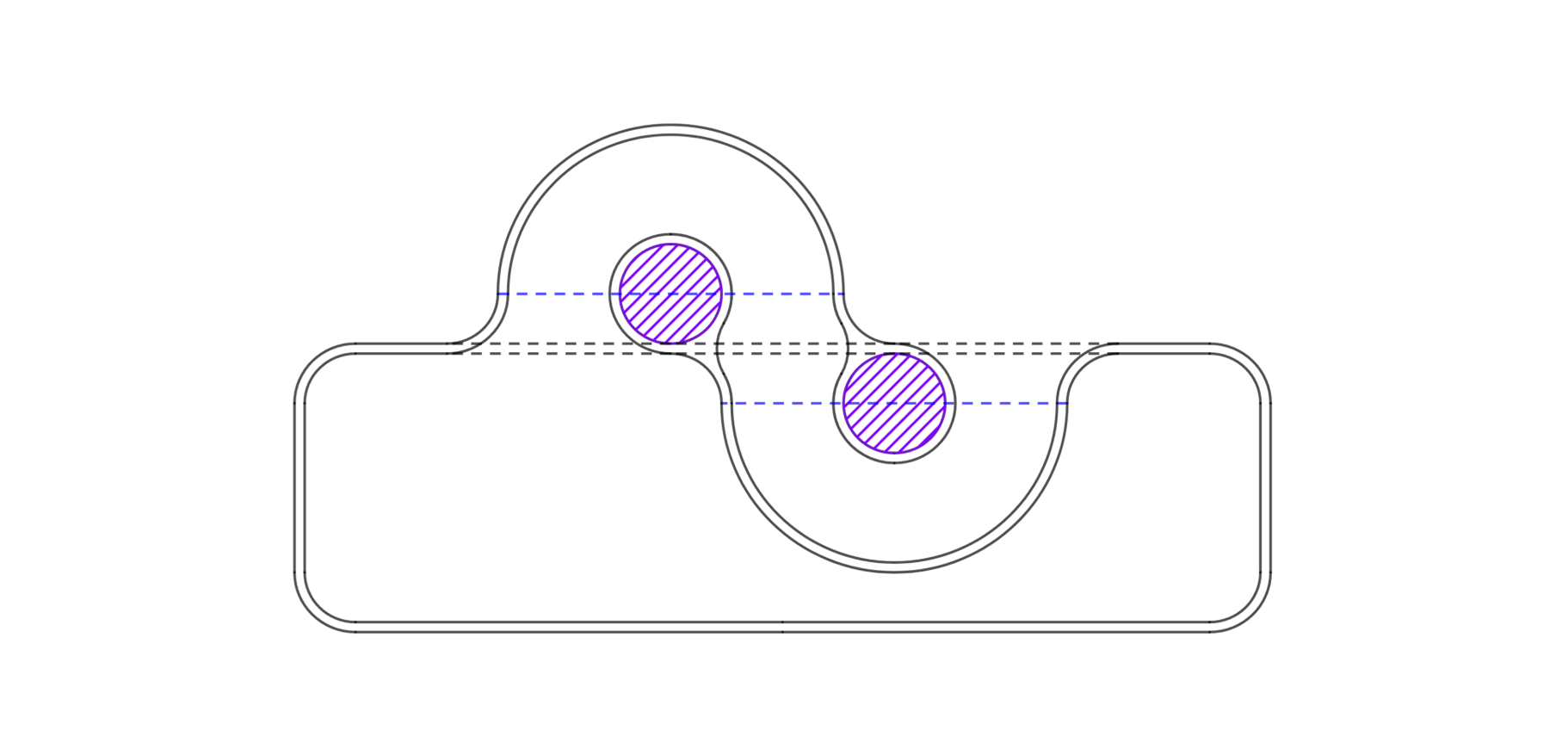}
}
\caption{The white part of the surface  is not important; the volume of the jerrycan near this part can be made arbitrarily small.
If we assume that thinness is zero, then
(1) red and blue parts are made from
pieces of tori obtained by rotating unit circles whose centers are at distances $2$ and $2-\sqrt{3}$ from the axis of rotation, respectively;
(2) yellow pieces are cut from the unit sphere;
(3) the black part is made from a round cylinder;
(4) green semidiscs are flat.
}
\label{figure counter example first appear}
\end{figure}

Our example (let us call it jerrycan) and its cross-section by the plane of its mirror symmetry are shown in Figure~\ref{figure counter example first appear}.
The jerrycan is a twisted version of Lagunov's fishbowl with several engineering improvements by Anton Petrunin and Rostislav Matveev;
my contribution is verifying that everything works.
This last step was also done independently by Matthew Bolan in a slightly different way; see \cite{anton2022mathoverflow}.

The annotated colors almost completely describe the shape.
Table~\ref{figure canonical} shows more details of the main part of the example; this is what remains after removing the several infinitely thin pieces.
In particular,
it shows the geometry of the S-tube that goes from the inside to the outside of the jerrycan.

\begin{longtable}{|c|l|}
\caption{The main part}
\label{figure canonical}\\ 
\hline
\begin{minipage}{40mm}
\vskip3mm
\centering
\includegraphics[width=0.9\linewidth]{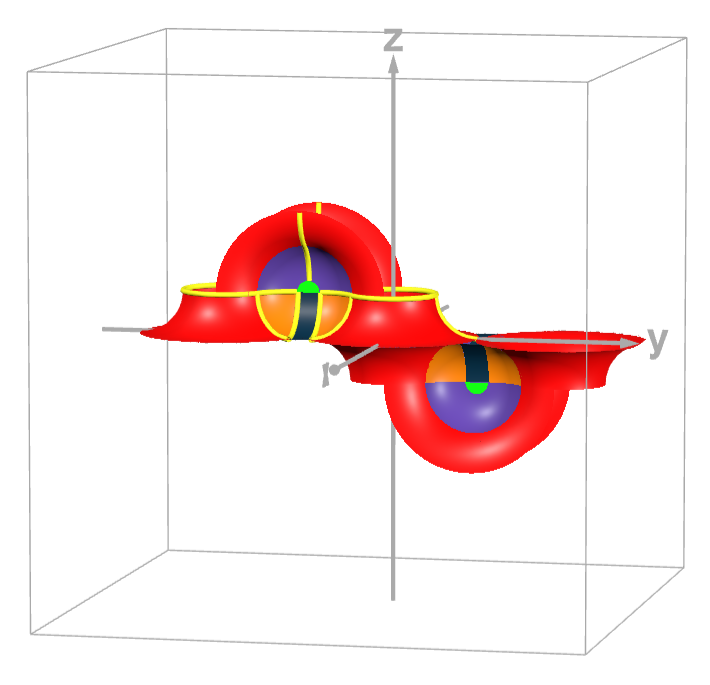}\
\\ \
\end{minipage}
&
\begin{minipage}{70mm}

\

The main part of the jerrycan contains all its volume, which is
\[2a+2b+c <\tfrac{4}{3}\pi;\]
the values $a$, $b$ and $c$ are defined below.

\

\end{minipage}
\\
\hline
\begin{minipage}{40mm}
\vskip3mm
\centering
\includegraphics[width=0.9\linewidth]{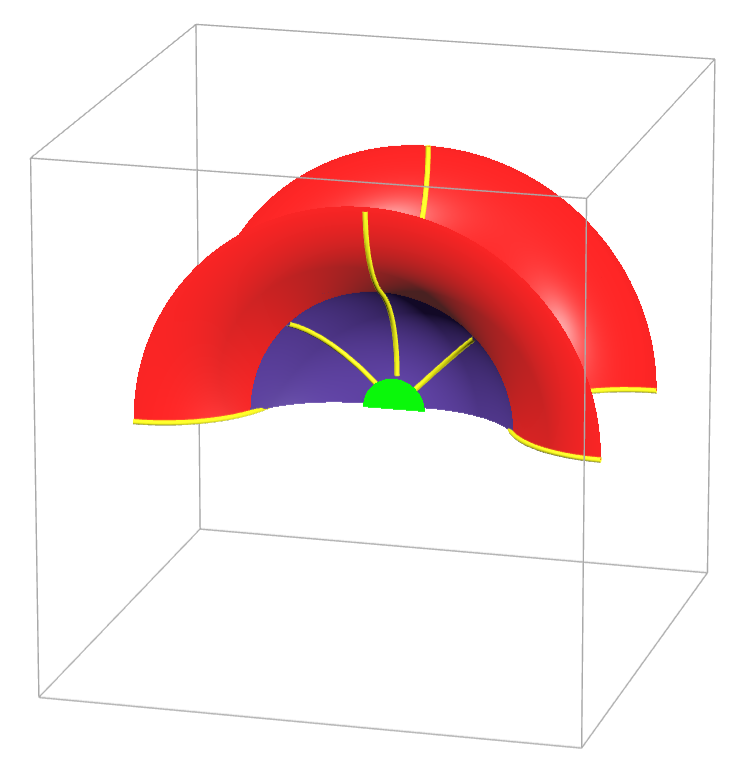}\
\\ \
\end{minipage}
&
\begin{minipage}{70mm}

\

The upper half of the cork. The main part has two such pieces; the volume of each is
\[a=\pi\cdot(2-\tfrac2{\sqrt3})\cdot(\sqrt3-\tfrac\pi2)< 0.5.\]

\

\end{minipage}
\\
\hline
\begin{minipage}{40mm}
\vskip3mm
\centering
\includegraphics[width=0.9\linewidth]{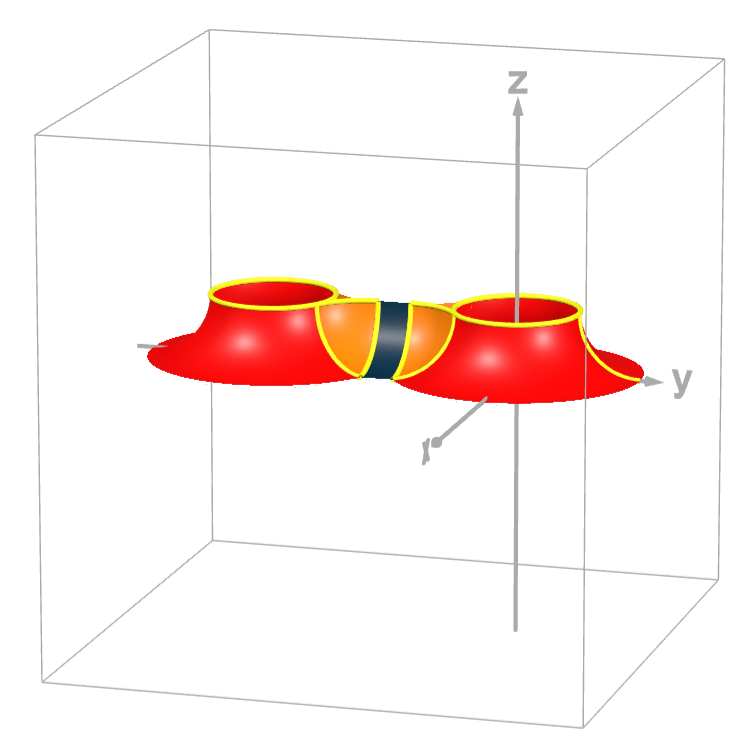}\
\\ \
\end{minipage}
&
\begin{minipage}{70mm}

\
The lower half of the cork, between the red toric slices. The main part has two such pieces; the volume of each is
\[b=8-2\sqrt{3}+\pi\sqrt{3}+\frac{\pi^2}{3}-4\pi<0.8.\]

\

\end{minipage}
\\
\hline
\begin{minipage}{40mm}
\vskip3mm
\centering
\includegraphics[width=0.9\linewidth]{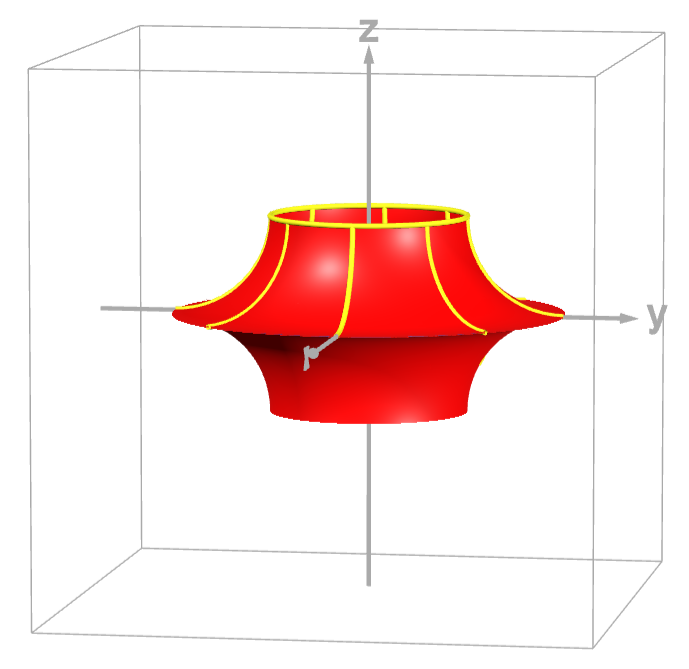}\
\\
\includegraphics[width=0.9\linewidth]{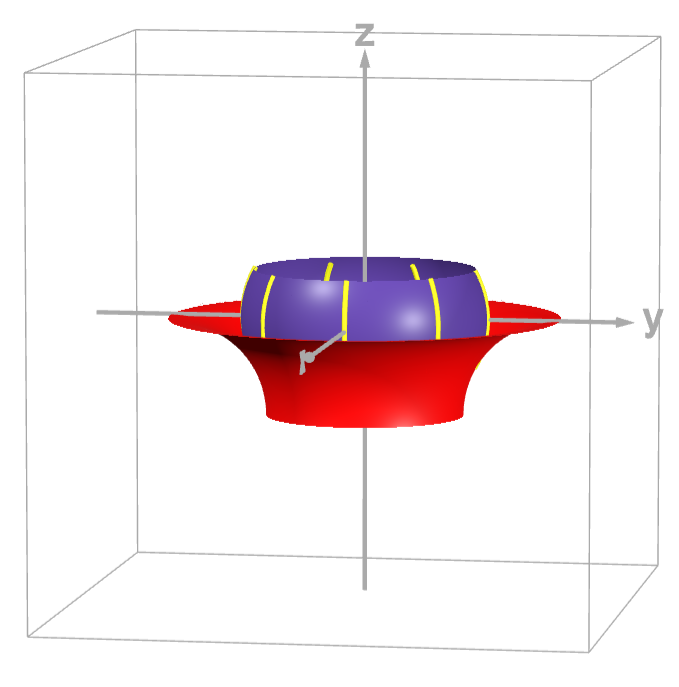}\
\\
\end{minipage}
&
\begin{minipage}{70mm}

\

The central piece and its tube; its volume is
\[c=
4\pi\sqrt{3}-2\pi-2\pi^2+\frac{\pi^2}{\sqrt{3}} <1.5.
\]

\

\end{minipage}
\\
\hline
\end{longtable}

The volume calculations require patience, but almost no ingenuity.

\begin{figure}[htbp]
\centering
\includegraphics[width=0.8\linewidth]{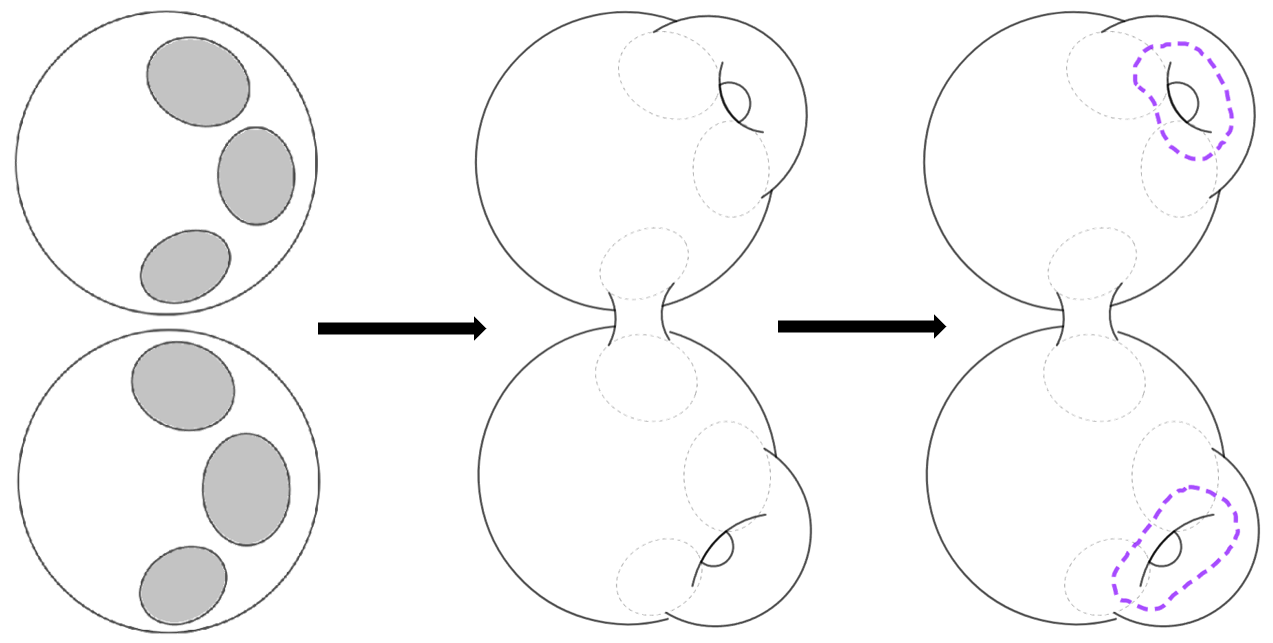}
\label{figure topology}
\end{figure}

One should also check that the surface is indeed a topological sphere. This is concisely illustrated in Figure~\ref{figure topology}. One starts with two spheres with three holes on each, add two handles and a tube, cut the resulting surface along the purple circles and glue four discs.

Finally, our jerrycan is only $C^{1,1}$; to make it smooth one has to inflate it a bit so that the curvature becomes strictly less than 1 and then apply smoothing (any smoothing procedure will do the trick).

\section{Remarks}

Our jerrycan does not minimize the volume.
Finding the optimal bound must be a very hard problem.
Moreover, even just guessing the optimal shape seems to be challenging.
However, we conjecture that $c=4\pi\sqrt{3}-2\pi-2\pi^2+\frac{\pi^2}{\sqrt{3}}$
gives the optimal value for closed surfaces of even positive genus (sphere with 2 handles, sphere with 4 handles, and so on).
The value $c$
is the volume of the tube shown in the last row of Table~\ref{figure canonical}.
Given an even genus $g>0$,
it is possible to construct a surface of genus $g$
with normal curvatures at most 1 and enclosed volume arbitrarily close to $c$ (an exercise), but it is unclear if one can achieve a smaller volume. For odd genus $g\ge 3$, it is possible to produce an example with volume arbitrarily close to $c+2a$, and this might be the optimal bound.
The torus case seems to be more complicated.
By the way, Lagunov claims \cite{lagunov1959russian}\cite{lagunovfet1961russianII} that in the torus case the maximal radius of the enclosed ball can be arbitrarily close to $\frac{2\sqrt{3}}{3}-1$.
This result was attributed to Valentin Diskant, but it was never published, and I have not been able to construct such an example myself.

My earlier result \cite{hongda2025ball} suggests the following question:

\smallskip

\textit{Let $V_{g,R}$ be the minimal enclosed volume among all $C^2$-smooth surfaces in $\RR^3$ of genus $g$ which have normal curvatures $\leq 1$ in absolute value and which are contained within a sphere of radius $R$. Then for every $g\geq0$, does there exist a sufficiently large $R_g$ such that $V_g = V_{g,R_g}$?}

\smallskip

\textbf{Acknowledgments:} I would like to thank Dmitri Burago and Anton Petrunin for their patience and help.

\bibliography{D:/latex/ref}
\bibliographystyle{plainurl}
\end{document}